\documentclass[12pt, a4]{amsart}

\usepackage{amscd, amsmath, amssymb}
\usepackage{fullpage}
\usepackage{mathrsfs,epsfig}
\usepackage{float}
\usepackage{algorithm}
\usepackage{algpseudocode}
\usepackage{matlab-prettifier}
\usepackage{subcaption}

\theoremstyle{plain}
\newtheorem{thm}{Theorem}[section]

\theoremstyle{definition}

\newtheorem{ex}[thm]{Example}

\numberwithin{equation}{section}

\newcommand{\R}{\mathbb{R}}
\newcommand{\N}{\mathbb{N}}

\newcommand{\norm}[1]{\|#1\|_\infty}
\DeclareMathOperator{\sgn}{sgn}

\usepackage{color}

\begin{document}
	
	\title[Eigenvalues of Neumann BVPs]{An eigenvalue result for Neumann BVPs with functional terms}
	
	\date{}
	
	\author[G. A. Veltri]{Giuseppe Antonio Veltri}
	\address{Giuseppe Antonio Veltri, Dipartimento di Scienze Matematiche e Informatiche, Scienze Fisiche e Scienze della Terra, Universit\`{a} degli Studi di Messina, 98100 Messina, Italy}%
	\email{giuseppe.veltri@studenti.unime.it}%
	
	\begin{abstract}
		We study the existence and localization of eigenvalue-eigenfunction pairs for parameter-dependent Neumann BVPs with a functional term. By reformulating the problems as a Hammerstein integral equation, we apply an existence and localization result and propose a convergent fixed-point iteration scheme. Finally, two pseudocodes and a MATLAB implementation are provided to numerically approximate the eigenvalues and validate the theoretical localization bounds. We also illustrate an approximation of the eigenfunctions for a fixed norm.
	\end{abstract}
	
	
	
	\maketitle
	
	\section{Introduction}
	In this paper, we aim to show the existence of couples $(\lambda,u)$ that solve the Neumann boundary value problem (BVP)
	\begin{equation}\label{Neu-BVP}
		\begin{cases}
			\epsilon u''(t)+\omega^2u(t)=\lambda f(t,u(t),H[u]),&t\in(0,1),\\
			u'(0)=u'(1)=0,
		\end{cases}
	\end{equation}
	where $\varepsilon=\pm 1$, $\omega$ is a suitable positive number, $f$ is a suitable function and $H$ a suitable functional. The existence of solutions for Neumann BVPs has been object of interest by a number of authors. For instance, Bonanno and Pizzimenti \cite{BonPiz12} studied the existence of positive solutions of
	\[
	\begin{cases}
		-u''(t)+u(t)=\lambda f(t,u(t)),&t\in(0,1),\\
		u'(0)=u'(1)=0,
	\end{cases}
	\]
	where $f$ is $L^1$-Caratheodory and $\lambda>0$. Their approach relies on variational methods. The same approach has been followed by Pizzimenti and Sciammetta \cite{PizSci10} in the searching of positive solutions of the problem
	\[
	\begin{cases}
		-u''(t)+uh(u'(t))=\lambda \alpha(t)f(u(t))h(u'(t)),&t \in(0,1),\\
		u'(0)=u'(1)=0,
	\end{cases}
	\]
	with $\lambda>0$, $\alpha\in C([0,1],(0,+\infty))$ and $f,h\in C(\R,\R)$. Other similar parameter-dependent Neumann BVPs have been studied in \cite{AfrHagSol23,AveBon04,BisBonRad11,BisBonRad13} and references therein.
	
	The non-parameter-dependent (i.e. $\lambda=1$) problem has been subject of interest too. For example, Henderson and Kosmatov \cite{HenKos15}, discussed the solvability of the BVP
	\[
	\begin{cases}
		-u''(t)+\omega^2u(t)=g(t,u(t)),&t \in(0,1),\\
		u'(0)=u'(1)=0,
	\end{cases}
	\]
	while Infante, Pietramala and Tojo \cite{InfPieToj16} discussed about the BVPs
	\[
	\begin{cases}
		\pm u''(t)+\omega^2u(t)=f(t,u(t)),&t\in(0,1),\\
		u'(0)=u'(1)=0.
	\end{cases}
	\]
	In both \cite{HenKos15} and \cite{InfPieToj16}, the authors rewrote the studied BVPs into a Hammerstein integral equation and their approach relies on topological fixed-point theory. Topological methods have been applied in the study of existence of solutions also by other authors in similar Neumann BVPs. We recall, for istance, the papers of Wang, Yu and Zhang~\cite{WanYuZha10}, Yao~\cite{Yao12}, Zhai and Zhang~\cite{ZhaZha10}, Nkashama and Santanilla~\cite{NkaSan90}, Zhilong~\cite{Zhi09}, and references therein.
	
	This manuscript is organized as follows: in Section 2, we recall a recent result by Infante and Veltri~\cite{InfVel25}, who applied it to the context of Dirichlet and mixed boundary conditions~(BCs), on the existence and localization of eigenvalues and eigenfunctions of Hammerstein integral equations of the form
	\begin{equation}\label{HIE}
		u(t) = \lambda\int_0^1 k(t,s) f(s,u(s),H[u])\,ds.
	\end{equation}
	Note that integral equations of a form similar to~\eqref{HIE} have been studied by a number of authors, we refer the reader to the Introduction of~\cite{InfVel25}.
	We will rewrite, in both the cases $\epsilon=\pm 1$, the BVP~\eqref{Neu-BVP} as a Hammerstein integral equation of the form \eqref{HIE} and apply the result to two examples to show the existence of eigenpairs of the equation.
	In Section~3, we develop a fixed-point iteration and prove that it converges, up to a subsequence, to a solution pair $(\lambda^*,u^*)$ by using the Bolzano-Weierstrass and the Arzelà-Ascoli theorems.
	Finally, in Section~4, we provide two pseudocodes and a MATLAB implementation of the proposed iteration that let us approximate the eigenvalues and plot them along with a theoretical localization of the eigenvalues that depends on the given norm of their associated eigenfunction. We also illustrate an approximation of the eigenfunctions for a fixed norm.
	
	\section{Main results}
	Here, we focus on the solvability of Neumann BVPs of the form
	\begin{equation}\label{Neu-BVP-2}
		\begin{cases}
			\epsilon u''(t)+\omega^2u(t) = \lambda f(t,u(t),H[u]), &t \in (0,1),\\
			u'(0)=u'(1)=0,
		\end{cases}
	\end{equation}
	where $\epsilon=\pm 1$, $\omega>0$ in the case $\epsilon=-1$ and $\omega \in (0,\pi/2]$ in the case $\epsilon=1$, $f$ is a suitable function and $H$ a suitable functional.
	It is known, see for example~\cite{InfPieToj16}, that we can rewrite the BVP~\eqref{Neu-BVP-2} into an equivalent Hammerstein integral equation of the form
	\begin{equation}\label{HIE2}
		u(t) = \lambda\int_0^1 k(t,s) f(s,u(s),H[u])\,ds,    
	\end{equation}
	where $k$ is the Green's function associated to \eqref{Neu-BVP-2}, that is
	\begin{equation}\label{greenneumannminus}
		k(t,s)=\frac{1}{\omega \sinh(\omega)}
		\begin{cases}
			\cosh(\omega(1-s))\cosh(\omega t), &0\leq t\leq s\leq1, \\
			\cosh(\omega(1-t))\cosh(\omega s), &0\leq s\leq t\leq1,
		\end{cases}
	\end{equation}
	in the case $\epsilon=-1$, and (with abuse of notation)
	\begin{equation}\label{greenneumannplus}
		k(t,s)=\frac{1}{\omega \sin(\omega)}
		\begin{cases}
			\cos(\omega(1-s))\cos(\omega t), &0\leq t\leq s\leq1, \\
			\cos(\omega(1-t))\cos(\omega s), &0\leq s\leq t\leq1.
		\end{cases}
	\end{equation}
	in the case $\epsilon=1$. In this latter case, we consider $\omega \in (0,\pi/2]$ because, for such values of $\omega$, the Green's function $k$ is always non-negative (as a reference the reader may take Lemma~5.1 in \cite{InfPieToj16}).
	We use the following notation:
	\[
	B_\rho := B_\rho(C([0,1])) := \{u \in C[0,1] \mid \norm{u} < \rho \} \quad\text{for every $\rho>0$}.
	\]
	
	The main theoretical tool that we will use in our reasoning it is the following and is based on a version of the classical Birkhoff-Kellog Theorem:
	\begin{thm}\label{HIEresult}\cite{InfVel25}
		Let \( \rho \in (0, +\infty) \) and assume that:
		\begin{enumerate}
			\item \( k : [0,1]^2 \rightarrow [0, +\infty) \) is continuous.
			\item \( H : \overline{B}_\rho \rightarrow \mathbb{R} \) is continuous and there exist \( \underline{H_\rho}, \overline{H_\rho} \in \mathbb{R} \) such that 
			\[ \underline{H_\rho} \leq H[u] \leq \overline{H_\rho}, \text{ for every } u \in \overline{B_\rho}. \]
			\item \( f : \Pi_\rho \subset \mathbb{R}^3 \rightarrow \mathbb{R} \) is continuous, where \(\Pi_\rho := [0,1] \times [-\rho,\rho] \times [\underline{H_\rho}, \overline{H_\rho}]\).
			\item There exist two continuous functions \( \underline{f_\rho}, \overline{f_\rho} : [0,1] \rightarrow \mathbb{R} \) such that
			\[ \underline{f_\rho}(t) \leq f(t,u,v) \leq \overline{f_\rho}(t), \text{ for every } (t,u,v) \in \Pi_\rho.\]
			\item[$(5)$] Let 
			\begin{align*}
				\underline{F_\rho}(t) := \int_0^1 k(t,s)\underline{f_\rho}(s)\,ds,\quad
				\overline{F_\rho}(t) := \int_0^1 k(t,s)\overline{f_\rho}(s)\,ds,
			\end{align*}
			and assume that at least one of the following conditions holds:
			\begin{itemize}
				\item[$(5a)$] There exists $t_\rho \in [0,1]$ such that $\overline{F_\rho}(t_\rho)<0$.
				\item[$(5b)$] There exists $t_\rho \in [0,1]$ such that $\underline{F_\rho}(t_\rho)>0$.
			\end{itemize}
		\end{enumerate}
		Then there exist $\lambda_\rho^+>0$, $\lambda_\rho^-<0$ and $u^+_\rho,u^-_\rho\in\partial B_\rho$ such that $(\lambda_\rho^+,u^+_\rho)$ and $(\lambda_\rho^-,u^-_\rho)$ solve the integral equation~\eqref{HIE}.
		
		Furthermore, assume that the couple $(\lambda_\rho,u_\rho)$ satisfies the integral equation~\eqref{HIE2}. Then the following is true.
		\begin{itemize}
			\item[$(6a)$] If $(5a)$ holds, then we have the estimate
			$|\lambda_\rho| \leq -\dfrac{\rho}{\overline{F_\rho}(t_\rho)}.$\vspace{1mm}
			
			\item[$(6b)$] If $(5b)$ holds, then we have the estimate $|\lambda_\rho| \leq \dfrac{\rho}{\underline{F_\rho}(t_\rho)}$.
		\end{itemize}
	\end{thm}
	
	We now show, in two examples, the applicability of Theorem~\ref{HIEresult}.
	\begin{ex}\label{Exampleminus}
		Consider the BVP
		\begin{equation}\label{hnpeq}
			\begin{cases}
				-u''(t) + u(t)= \lambda  \dfrac{\sin\left(\frac{3}{2}\pi t\right)e^{u(t)}}{\int_0^1 e^{u(x)} \, dx}, & t\in(0,1),\\
				u'(0)=u'(1)=0.
			\end{cases}
		\end{equation}
		To the BVP~\eqref{hnpeq} we associate the Hammerstein integral equation
		\begin{equation}\label{HIEneumann}
			u(t) = \lambda \int_{0}^{1} k(t,s) \dfrac{\sin\left(\frac{3}{2}\pi s\right)e^{u(s)}}{\int_0^1 e^{u(x)} \, dx} \, ds,
		\end{equation}
		where $k(t,s)$ is as in \eqref{greenneumannminus}, that is
		\begin{equation}
			k(t,s)=\frac{1}{\sinh(1)}
			\begin{cases}
				\cosh(1-s)\cosh(t), &0\leq t\leq s\leq1, \\
				\cosh(1-t)\cosh(s), &0\leq s\leq t\leq1.
			\end{cases}
		\end{equation}
		Note that $\sin\left(\frac{3}{2}\pi t\right)$ is non-negative in $\left[ 0, 2/3 \right]$ and non-positive in $\left[ 2/3,1 \right]$, so we may consider
		\begin{align}
			\underline{F_\rho}(t) &= e^{-2\rho} \int_0^{\frac{2}{3}} k(t,s) \sin\left(\frac{3}{2}\pi s\right) \, ds + e^{2\rho}\int_{\frac{2}{3}}^1 k(t,s) \sin\left(\dfrac{3}{2}\pi s\right) \, ds,\\
			\overline{F_\rho}(t) &= e^{2\rho} \int_0^{\frac{2}{3}} k(t,s) \sin\left(\frac{3}{2}\pi s\right) \, ds + e^{-2\rho}\int_{\frac{2}{3}}^1 k(t,s) \sin\left(\dfrac{3}{2}\pi s\right) \, ds.
		\end{align}
		A direct calculation yields
		\[
		\underline{F_\rho}(t) =
		\begin{cases}
			e^{-2\rho} \mathcal{A}(t) + e^{2\rho} \mathcal{B}(t) , & \text{if}\quad 0 \leq t \leq \dfrac{2}{3},\vspace{1.5mm}\\
			e^{-2\rho} \mathcal{C}(t)  + e^{2\rho}\mathcal{D}(t), & \text{if}\quad \dfrac{2}{3} \leq t \leq 1,
		\end{cases}
		\]
		and
		\[
		\overline{F_\rho}(t) =
		\begin{cases}
			e^{2\rho} \mathcal{A}(t) + e^{-2\rho} \mathcal{B}(t) , & \text{if}\quad 0 \leq t \leq \dfrac{2}{3},\vspace{1.5mm}\\
			e^{2\rho} \mathcal{C}(t) + e^{-2\rho} \mathcal{D}(t), & \text{if}\quad \dfrac{2}{3} \leq t \leq 1,
		\end{cases}
		\]
		where, if we set $b=3\pi/2$ and $q=1/[(1+b^2)\sinh(1)]$,
		\begin{align*}
			\mathcal{A}(t) &:= q \sin(bt) \sinh(1) + qb \left( \cosh(1-t) + \cosh\left(\frac{1}{3}\right) \cosh(t) \right),\\
			\mathcal{B}(t) &:= -qb \cosh\left(\frac{1}{3}\right)\cosh(t),\\
			\mathcal{C}(t) &:= qb \left( \cosh\left(\frac{2}{3}\right) + 1 \right) \cosh(1-t),\\
			\mathcal{D}(t) &:= q \left( \sinh(1)\sin(bt) - b \cosh\left(\frac{2}{3}\right)\cosh(1-t) \right).
		\end{align*}
		Note that, for every fixed $\rho>0$, the functions $\underline{F_\rho}$ and $\overline{F_\rho}$ are $C^1$ and satisfy the homogeneous Neumann BCs. Figure~\ref{F_low-F_upp-Neumann-minus} illustrates the graphs of the functions $\underline{F_\rho}(t)$ and $\overline{F_\rho}(t)$ in the cases $\rho=0.1$, $\rho=0.15$ and $\rho=1$.
		\begin{figure}[H]
			\centering
			\makebox[\textwidth][c]{\includegraphics[width=1.1\linewidth]{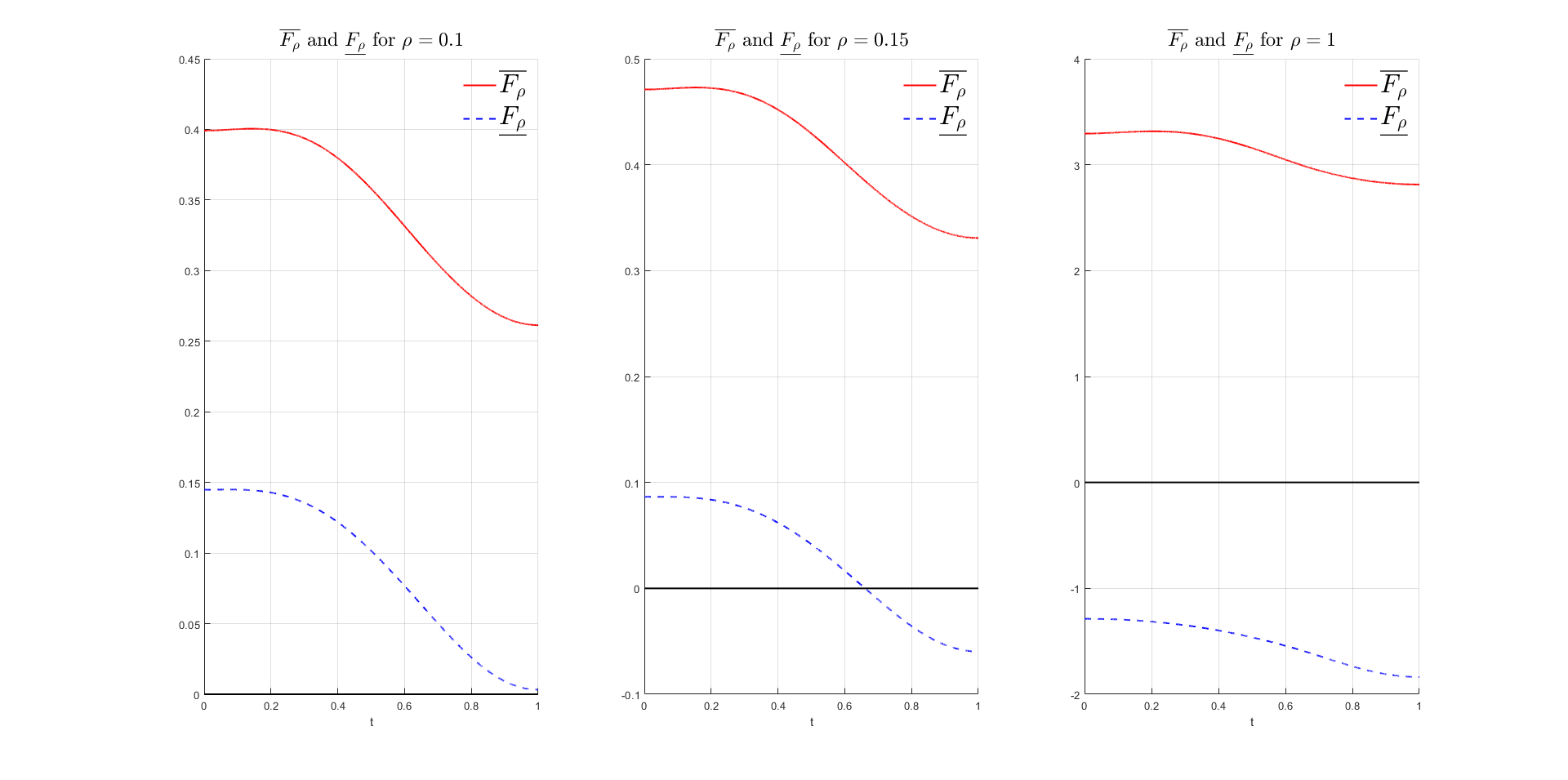}}
			\caption{Plots of $\underline{F_\rho}$ and $\overline{F_\rho}$ in the cases $\rho=0.1$, $\rho=0.15$ and $\rho=1$, respectively.}\label{F_low-F_upp-Neumann-minus}
		\end{figure}
		One can show that $\overline{F_\rho}$ is a positive function independently of $\rho$, and this means that, to apply Theorem~\ref{HIEresult}, we have to focus on the function $\underline{F_\rho}$.\\
		In particular, we have to search for the values of $\rho$ such that the function $\underline{F_\rho}$ achieves a positive maximum, that is an equivalent condition to $(5b)$ in Theorem~\ref{HIEresult}, and to find the absolute maximum point $t_\rho$ of $\underline{F_\rho}$, so that we can evaluate the best estimate possible in $(6b)$.
		
		First of all, note that $\mathcal{B}(x)<0$ and $\mathcal{D}(y)<0$ for every $x \in [0,2/3]$ and $y \in [2/3,1]$.
		
		Fix $\rho>0$. Assume that there exists $t_\rho\in[0,2/3]$ such that $\underline{F_\rho}(t_\rho)>0$, that is
		\[
		e^{-2\rho}\mathcal{A}(t_\rho)+e^{2\rho}\mathcal{B}(t_\rho)>0,
		\]
		or, equivalently,
		\[
		e^{4\rho} < -\mathcal{A}(t_\rho)/\mathcal{B}(t_\rho).
		\]
		This inequality yields
		\[
		\rho < \frac{1}{4} \log\left(-\frac{\mathcal{A}(t_\rho)}{\mathcal{B}(t_\rho)}\right),
		\]
		which is equivalent to say that
		\[
		\rho < \frac{1}{4} \max_{[0,2/3]}\log\left(-\frac{\mathcal{A}}{\mathcal{B}}\right) = \frac{1}{4} \log\left(\max_{[0,2/3]}\left(-\frac{\mathcal{A}}{\mathcal{B}}\right)\right)=:\rho_1.
		\]
		On the other hand, assume that there is $t_\rho\in[2/3,1]$ such that $\underline{F_\rho}(t_\rho)>0$, i.e.
		\[
		e^{-2\rho}\mathcal{C}(t_\rho)+e^{2\rho}\mathcal{D}(t_\rho)>0,
		\]
		that is
		\[
		e^{4\rho} < -\mathcal{C}(t_\rho)/\mathcal{D}(t_\rho).
		\]
		This inequality implies that
		\[
		\rho < \frac{1}{4} \log\left(-\frac{\mathcal{C}(t_\rho)}{\mathcal{D}(t_\rho)}\right),
		\]
		that is equivalent to the inequality
		\[
		\rho < \frac{1}{4} \max_{[2/3,1]}\log\left(-\frac{\mathcal{C}}{\mathcal{D}}\right) = \frac{1}{4} \log\left(\max_{[2/3,1]}\left(-\frac{\mathcal{C}}{\mathcal{D}}\right)\right) =: \rho_2.
		\]
		Therefore, $\underline{F_\rho}$ attains a positive maximum if and only if $0<\rho<\max\{\rho_1,\rho_2\}=:\rho_0$. Since the equations
		\[
		\left(-\frac{\mathcal{A}}{\mathcal{B}}\right)'(t)=0\quad\text{and}\quad\left(-\frac{\mathcal{C}}{\mathcal{D}}\right)'(t)=0
		\]
		are transcendental, we have to find $\rho_1$ and $\rho_2$, and so $\rho_0$, numerically.
		
		Fix $\rho\in(0,\rho_0)$. Assume that there exists an absolute maximum point $t_\rho \in (0,1)$ of $\underline{F_\rho}$. This implies that $\underline{F_\rho}'(t_\rho)=0$, where
		\[
		\underline{F_\rho}'(t_\rho) =
		\begin{cases}
			e^{-2\rho} \mathcal{A}'(t_\rho) + e^{2\rho} \mathcal{B}'(t_\rho) , & \text{if}\quad 0 < t_\rho \leq \dfrac{2}{3},\vspace{1.5mm}\\
			e^{-2\rho} \mathcal{C}'(t_\rho)  + e^{2\rho}\mathcal{D}'(t_\rho), & \text{if}\quad \dfrac{2}{3} \leq t_\rho < 1.
		\end{cases}
		\]
		Hence, we get two equations to solve:
		\[
		\begin{cases}
			\mathcal{A}'(t_\rho) + e^{4\rho} \mathcal{B}'(t_\rho) = 0, & \text{if}\quad 0 < t_\rho \leq \dfrac{2}{3},\vspace{1.5mm}\\
			\mathcal{C}'(t_\rho)  + e^{4\rho}\mathcal{D}'(t_\rho) = 0, & \text{if}\quad \dfrac{2}{3} \leq t_\rho < 1.
		\end{cases}
		\]
		Both these equations are transcendental, thus we cannot solve them analytically. This implies that we cannot compute $\max_{[0,1]} \underline{F_\rho}$ unless we do it numerically. This means that we have to compute both $\rho_0$ and the estimate in $(6b)$ in Theorem~\ref{HIEresult} with numerical methods. Figure~\ref{eigen_approx_minus} in the last section of this manuscript represents the plot of the approximation of both the positive and negative eigenvalues together with the estimates in $(6b)$. The vertical line $\rho=\rho_0$ is also present.
	\end{ex}
	
	\begin{ex}\label{Exampleplus}
		Consider the BVP
		\begin{equation}\label{hnpeqplus}
			\begin{cases}
				u''(t) + \dfrac{\pi^2}{4}u(t)= \lambda  \dfrac{\sin\left(\frac{3}{2}\pi t\right)e^{u(t)}}{\int_0^1 e^{u(x)} \, dx}, & t\in(0,1),\\
				u'(0)=u'(1)=0.
			\end{cases}
		\end{equation}
		To the BVP~\eqref{hnpeq} we associate the Hammerstein integral equation
		\begin{equation}
			u(t) = \lambda \int_{0}^{1} k(t,s) \dfrac{\sin\left(\frac{3}{2}\pi s\right)e^{u(s)}}{\int_0^1 e^{u(x)} \, dx} \, ds,
		\end{equation}
		where $k(t,s)$ is as in \eqref{greenneumannplus}, that is
		\begin{equation}
			k(t,s)=\frac{2}{\pi}
			\begin{cases}
				\cos\left(\dfrac{\pi}{2}(1-s)\right)\cos\left(\dfrac{\pi}{2}t\right), &0\leq t\leq s\leq1, \\
				\cos\left(\dfrac{\pi}{2}(1-t)\right)\cos\left(\dfrac{\pi}{2}s\right), &0\leq s\leq t\leq1.
			\end{cases}
		\end{equation}
		As in Example~\ref{Exampleminus}, we may take
		\begin{align}
			\underline{F_\rho}(t) &= e^{-2\rho} \int_0^{\frac{2}{3}} k(t,s) \sin\left(\frac{3}{2}\pi s\right) \, ds + e^{2\rho}\int_{\frac{2}{3}}^1 k(t,s) \sin\left(\dfrac{3}{2}\pi s\right) \, ds,\\
			\overline{F_\rho}(t) &= e^{2\rho} \int_0^{\frac{2}{3}} k(t,s) \sin\left(\frac{3}{2}\pi s\right) \, ds + e^{-2\rho}\int_{\frac{2}{3}}^1 k(t,s) \sin\left(\dfrac{3}{2}\pi s\right) \, ds.
		\end{align}
		A direct calculation yields
		\[
		\underline{F_\rho}(t) =
		\begin{cases}
			e^{-2\rho} \mathcal{A}(t) + e^{2\rho} \mathcal{B}(t) , & \text{if}\quad 0 \leq t \leq \dfrac{2}{3},\vspace{1.5mm}\\
			e^{-2\rho} \mathcal{C}(t)  + e^{2\rho}\mathcal{D}(t), & \text{if}\quad \dfrac{2}{3} \leq t \leq 1,
		\end{cases}
		\]
		and
		\[
		\overline{F_\rho}(t) =
		\begin{cases}
			e^{2\rho} \mathcal{A}(t) + e^{-2\rho} \mathcal{B}(t) , & \text{if}\quad 0 \leq t \leq \dfrac{2}{3},\vspace{1.5mm}\\
			e^{2\rho} \mathcal{C}(t) + e^{-2\rho} \mathcal{D}(t), & \text{if}\quad \dfrac{2}{3} \leq t \leq 1,
		\end{cases}
		\]
		where, by setting $\omega=\pi/2$, $b=3\pi/2=3\omega$ and $q=1/(\omega\sin\omega)=1/\omega$,
		\begin{align*}
			\mathcal{A}(t) &:= -\frac{q}{b^2 - \omega^2} \left[ \frac{\sin(bt)}{q} - b \cos\left(\frac{\omega}{3}\right) \cos(\omega t) - b \cos(\omega - \omega t) \right],\\
			\mathcal{B}(t) &:= -\frac{qb}{b^2 - \omega^2} \cos\left(\frac{\omega}{3}\right) \cos(\omega t),\\
			\mathcal{C}(t) &:= \frac{qb}{b^2 - \omega^2} \left[ \cos\left(\frac{2\omega}{3}\right) + 1 \right] \cos(\omega - \omega t),\\
			\mathcal{D}(t) &:= -\frac{q}{b^2 - \omega^2} \left[ \frac{\sin(bt)}{q} + b \cos\left(\frac{2\omega}{3}\right) \cos(\omega - \omega t) \right],
		\end{align*}
		that is, by using the triple-angle formulas,
		\begin{align*}
			\mathcal{A}(t) 
			&= \frac{1}{4\pi^2} \left[ 8\sin^3\left(\frac{\pi}{2}t\right) + 3\sqrt{3} \cos\left(\frac{\pi}{2}t\right)\right]\\
			\mathcal{B}(t)
			&= -\frac{3\sqrt{3}}{4\pi^2} \cos\left(\frac{\pi}{2}t\right),\\
			\mathcal{C}(t)
			&= \frac{9}{4\pi^2} \sin\left(\frac{\pi}{2}t\right),\\
			\mathcal{D}(t)
			&= -\frac{1}{4\pi^2}\sin\left(\frac{\pi}{2}t\right) \left[ 9 - 8 \sin^2\left(\frac{\pi}{2}t\right) \right].
		\end{align*}
		As in the previous example, for every fixed $\rho>0$, the functions $\underline{F_\rho}$ and $\overline{F_\rho}$ are $C^1$ and satisfy the homogeneous Neumann BCs. Figure~\ref{F_upp-F_low-Neumann-plus} illustrates the graphs of the functions $\underline{F_\rho}(t)$ and $\overline{F_\rho}(t)$ in the cases $\rho=0.1$, $\rho=0.4$ and $\rho=1$.
		\begin{figure}[h]
			\centering
			\makebox[\textwidth][c]{\includegraphics[width=1.1\linewidth]{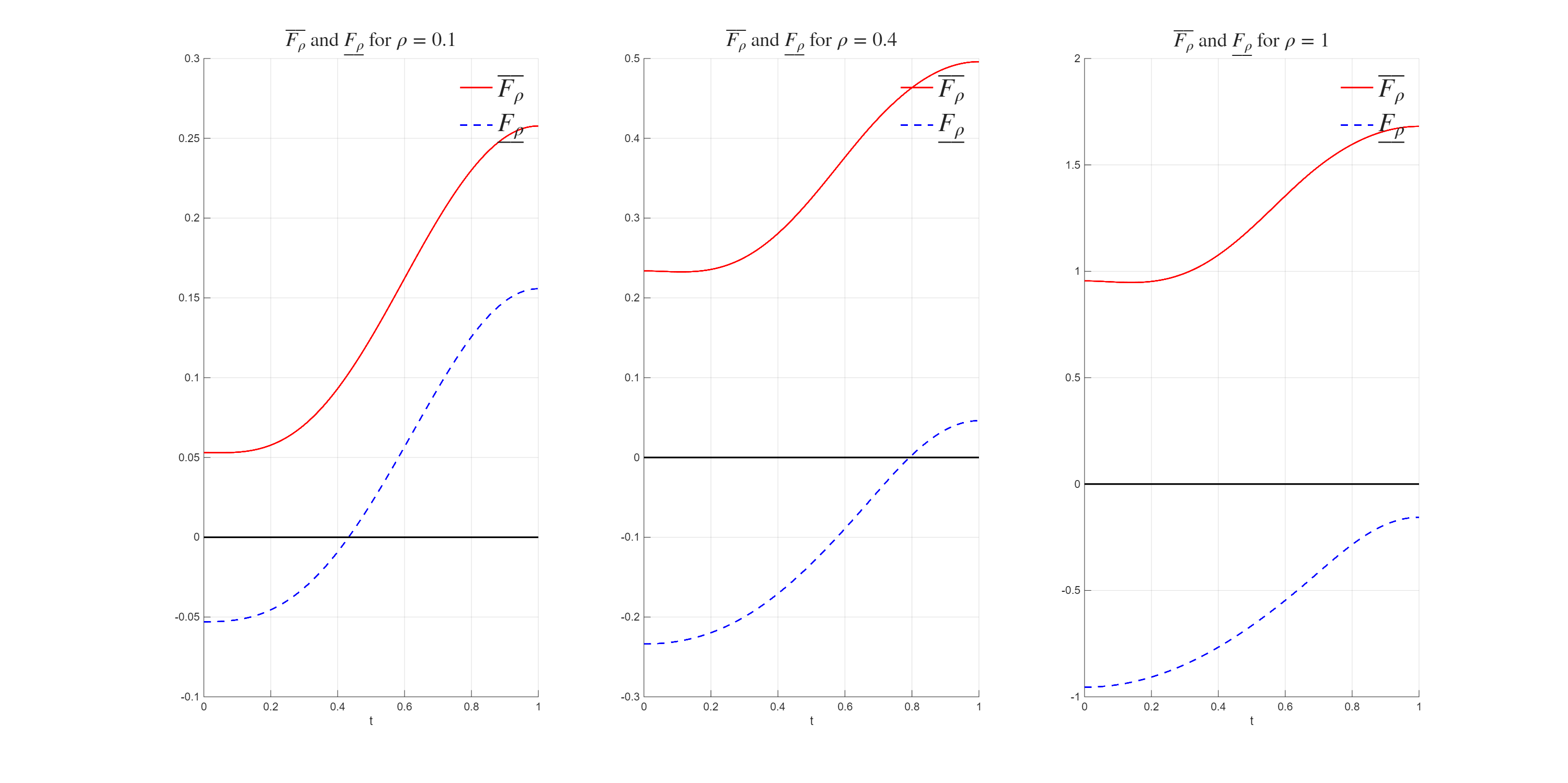}}
			\caption{Plots of $\underline{F_\rho}$ and $\overline{F_\rho}$ in the cases $\rho=0.1$, $\rho=0.4$ and $\rho=1$, respectively.}\label{F_upp-F_low-Neumann-plus}
		\end{figure}
		
		Note that $\overline{F_\rho}$ is a positive function independently of $\rho$, hence, to apply Theorem~\ref{HIEresult}, we have to concentrate on the function $\underline{F_\rho}$.
		
		Firstly, observe that $B(x),D(y)<0$ for every $x\in[0,2/3]$ and $y\in[2/3,1]$.
		Assume that $\rho>0$ is such that there exists $t_\rho \in [0,2/3]$ satisfying $\underline{F_\rho}(t_\rho)>0$, that is
		\[
		e^{4\rho} \mathcal{B}(t_\rho) > -\mathcal{A}(t_\rho),
		\]
		i.e.
		\[
		e^{4\rho} < -\frac{\mathcal{A}(t_\rho)}{\mathcal{B}(t_\rho)}.
		\]
		This latter inequality is equivalent to say that
		\[
		\rho < \frac{1}{4} \log\left(\max_{[0,2/3]}\left(-\frac{\mathcal{A}}{\mathcal{B}}\right)\right)=:\rho_1.
		\]
		Now, suppose that $\rho>0$ is such that there is $t_\rho \in [2/3,1]$ that satisfies $\underline{F_\rho}(t_\rho)>0$, that is
		\[
		e^{4\rho} \mathcal{D}(t_\rho) > -\mathcal{C}(t_\rho),
		\]
		that is equivalent to say that
		\[
		e^{4\rho} < -\frac{\mathcal{C}(t_\rho)}{\mathcal{D}(t_\rho)}.
		\]
		This inequality is equivalent to the following:
		\[
		\rho < \frac{1}{4} \log\left(\max_{[0,2/3]}\left(-\frac{\mathcal{C}}{\mathcal{D}}\right)\right)=:\rho_2.
		\]
		Therefore, $\underline{F_\rho}$ attains a positive maximum if and only if $0<\rho<\max\{\rho_1,\rho_2\}=:\rho_0$.
		
		Note that
		\begin{align*}
			\left(-\frac{\mathcal{A}(t)}{\mathcal{B}(t)}\right)'&= \left(\frac{8}{3\sqrt{3}}\tan\left(\frac{\pi}{2}t\right)\sin^2\left(\frac{\pi}{2}t\right)+1\right)'\\
			&=\frac{4\pi}{3\sqrt{3}}\left[\sin^2\left(\frac{\pi}{2}t\right)+\tan^2\left(\frac{\pi}{2}t\right)\right] \geq 0\\
		\end{align*}
		for every $t \in (0,2/3)$. This means that the maximum of $-\mathcal{A}/\mathcal{B}$ is
		\[
		-\frac{\mathcal{A}(2/3)}{\mathcal{B}(2/3)} = \left(\frac{8}{3\sqrt{3}}\tan\left(\frac{\pi}{3}\right)\sin^2\left(\frac{\pi}{3}\right)+1\right)=3,
		\]
		and this yields
		$$\rho_1=\log(\sqrt[4]{3}).$$
		On the other side,
		\[
		-\frac{\mathcal{C}(t)}{\mathcal{D}(t)}= \frac{9}{9-8\sin^2\left(\frac{\pi}{2}t\right)}
		\]
		attains its maximum when $9-8\sin^2\left(\frac{\pi}{2}t\right)$ attains its minimum, that is in $t=1$. Therefore, we get that
		\[
		\rho_2 = \log\sqrt[4]{\frac{9}{9-8\sin^2\left(\frac{\pi}{2}\right)}}=\log\sqrt{3}.
		\]
		Therefore, we have that $\rho_0=\log\sqrt{3}\approx0.5493$.
		
		Now, assume that $\rho \in (0,\rho_0)$. We want to compute $\max_{[0,1]}\underline{F_\rho}$. To do so, we need to find the maximum point $t_\rho$ of $\underline{F_\rho}$.
		
		Let us recall that
		\[
		\underline{F_\rho}'(t) =
		\begin{cases}
			e^{-2\rho} \mathcal{A}'(t) + e^{2\rho} \mathcal{B}'(t) , & \text{if}\quad 0 < t \leq \dfrac{2}{3},\vspace{1.5mm}\\
			e^{-2\rho} \mathcal{C}'(t)  + e^{2\rho}\mathcal{D}'(t), & \text{if}\quad \dfrac{2}{3} \leq t < 1,
		\end{cases}
		\]
		where
		\begin{align*}
			\mathcal{A}'(t) &= \frac{3}{\pi} \sin^2\left(\frac{\pi}{2}t\right)\cos\left(\frac{\pi}{2}t\right) - \frac{3\sqrt{3}}{8\pi} \sin\left(\frac{\pi}{2}t\right) \\
			\mathcal{B}'(t) &= \frac{3\sqrt{3}}{8\pi} \sin\left(\frac{\pi}{2}t\right) \\
			\mathcal{C}'(t) &= \frac{9}{8\pi} \cos\left(\frac{\pi}{2}t\right) \\
			\mathcal{D}'(t) &= \frac{3}{\pi} \sin^2\left(\frac{\pi}{2}t\right)\cos\left(\frac{\pi}{2}t\right) - \frac{9}{8\pi} \cos\left(\frac{\pi}{2}t\right)
		\end{align*}
		Since $t_\rho$ is the maximum point of $\underline{F_\rho}$, we have that $\underline{F_\rho}'(t_\rho)=0$, that is
		\begin{equation}\label{AuxiliaryEq}
			\begin{cases}
				\mathcal{A}'(t_\rho) + e^{4\rho} \mathcal{B}'(t_\rho) = 0, & \text{if}\quad 0 < t_\rho \leq \dfrac{2}{3},\vspace{1.5mm}\\
				\mathcal{C}'(t_\rho)  + e^{4\rho}\mathcal{D}'(t_\rho) = 0, & \text{if}\quad \dfrac{2}{3} \leq t_\rho < 1.
			\end{cases}
		\end{equation}
		Note that $\mathcal{B}'\neq 0$ in $(0,2/3]$ and $\mathcal{C}'\neq 0$ in $[2/3,1)$.
		In particular, $\eqref{AuxiliaryEq}$ reads as
		\[
		\begin{cases}
			e^{4\rho} = -\mathcal{A}'(t_\rho)/\mathcal{B}'(t_\rho), & \text{if}\quad 0 < t_\rho \leq \dfrac{2}{3},\vspace{1.5mm}\\
			e^{4\rho} = -\mathcal{C}'(t_\rho)/\mathcal{D}'(t_\rho), & \text{if}\quad \dfrac{2}{3} \leq t_\rho < 1.
		\end{cases}
		\]
		This implies, since $\rho>0$, that $-\mathcal{A}'(t_\rho)/\mathcal{B}'(t_\rho)>1$ or $-\mathcal{C}'(t_\rho)/\mathcal{D}'(t_\rho)>1$, depending on which sub-interval of $(0,1)$ $t_\rho$ belongs to, but these two inequalities are never attained. Indeed,
		\begin{align*}
			-\frac{\mathcal{A}'(t_\rho)}{\mathcal{B}'(t_\rho)} &= -\frac{\frac{3}{\pi} \sin^2\left(\frac{\pi}{2}t_\rho\right)\cos\left(\frac{\pi}{2}t_\rho\right) - \frac{3\sqrt{3}}{8\pi} \sin\left(\frac{\pi}{2}t_\rho\right)}{\frac{3\sqrt{3}}{8\pi} \sin\left(\frac{\pi}{2}t_\rho\right)} \\
			&= 1 - \frac{8}{\sqrt{3}} \sin\left(\frac{\pi}{2}t_\rho\right)\cos\left(\frac{\pi}{2}t_\rho\right) \\
			&= 1 - \frac{4}{\sqrt{3}} \sin(\pi t_\rho) < 1,
		\end{align*}
		if $t_\rho \in (0,2/3]$, and
		\begin{align*}
			-\frac{\mathcal{C}'(t_\rho)}{\mathcal{D}'(t_\rho)} &= -\frac{\frac{9}{8\pi} \cos\left(\frac{\pi}{2}t_\rho\right)}{\frac{3}{\pi} \sin^2\left(\frac{\pi}{2}t_\rho\right)\cos\left(\frac{\pi}{2}t_\rho\right) - \frac{9}{8\pi} \cos\left(\frac{\pi}{2}t_\rho\right)} \\
			&= -\frac{\frac{9}{8\pi} \cos\left(\frac{\pi}{2}t_\rho\right)}{\frac{9}{8\pi} \cos\left(\frac{\pi}{2}t_\rho\right) \left[ \frac{8}{3} \sin^2\left(\frac{\pi}{2}t_\rho\right) - 1 \right]} \\
			&= \frac{3}{3 - 8\sin^2\left(\frac{\pi}{2}t_\rho\right)} < 0 < 1,
		\end{align*}
		if $t_\rho \in [2/3,1)$.
		This means that the maximum point of $\underline{F_\rho}$ must be $t_\rho=0$ or $t_\rho=1$, and a direct computation shows that
		\[
		\max_{[0,1]}\underline{F_\rho} = \underline{F_\rho}(1)=\frac{9e^{-2\rho}-e^{2\rho}}{4\pi^2} > 0 \quad\text{for every }\rho\in(0,\rho_0).
		\]
	\end{ex}
	
	\section{A fixed-point iteration}
	In this Section, we will construct a fixed-point iteration that converges, up to a subsequence, to the couple $(u^*,\lambda^*)$ that solves the Hammerstein integral equation \eqref{HIE} in the case of Example~\ref{Exampleminus}, for every fixed $\rho \in (0,\rho_0)$, that is such that Theorem~\ref{HIEresult} holds. This construction follows the setting developed, with the help of Professor Giovanni Mascali, in the Master's degree thesis \cite{Thesis}. The novelty here is the tailoring to the case of Neumann BCs.
	
	Let us fix $\rho\in(0,\rho_0)$. We seek a couple $(u^*,\lambda^*) \in C[0,1]\times\R$ such that
	\[
	\begin{cases}
		u^*(t) = \lambda^* Tu^*(t) &\text{for every } t \in [0,1],\\
		\norm{u^*} = \rho.
	\end{cases}
	\]
	Consider an initial guess $(u_0,\lambda_0)$, where $\norm{u_0}=\rho$ and $\sgn(\lambda_0)=\sgn(\lambda^*)$, and define for every $n \in \N \cup \{0\}=: \N_0$
	\begin{gather*}
		\lambda_{n+1} = \sgn(\lambda^*)\rho / \|Tu_n\|_\infty,\\
		u_{n+1}(t) = \lambda_{n+1} Tu_n(t) \quad \text{for every } t \in [0,1].
	\end{gather*}
	Note that $\{u_n\}_n$ is uniformly bounded in $C[0,1]$ and $\{\lambda_n\}_n$ is bounded in $\R$. Indeed, we have that
	\[
	\norm{u_{n+1}} = |\lambda_{n+1}| \norm{Tu_n} = \rho \quad\text{for each }n\in\N_0
	\]
	and this implies that
	\[
	\norm{Tu_n} \geq \max_{[0,1]}\underline{F_\rho} > 0 \quad\text{for every }n\in\N_0.
	\]
	Therefore, $\{\lambda_n\}_n$ is well-defined and
	\[
	|\lambda_n| =
	\begin{cases}
		|\lambda_0|,&\text{if } n=0,\\
		\rho / \|Tu_{n-1}\|_\infty,&\text{if } n\geq 1.
	\end{cases}
	\]
	Therefore, for every $n \in \N_0$,
	\[
	|\lambda_n| \leq
	\max\Big\{|\lambda_0|,\rho/\max_{[0,1]}\underline{F_\rho}\Big\}=:c_\rho.
	\]
	Moreover, $\{u_n\}_n$ is an equicontinuous sequence of functions. Indeed, since $k$ is uniformly continuous, for every $\varepsilon>0$ there exists $\delta_1>0$ such that for each $t_1,t_2 \in [0,1]$ satisfying $|t_1-t_2|<\delta_1$ we have that $|k(t_1,s)-k(t_2,s)|<\varepsilon$ for every $s \in [0,1]$. Also, the continuity of $f$ in $\Pi_\rho$ yields the existence of $\max_{\Pi_\rho}|f|$. Therefore, we get that, for every $t_1,t_2 \in [0,1]$ such that $|t_1-t_2|<\delta_1$,
	\begin{align*}
		|u_{n+1}(t_1) - u_{n+1}(t_2)| & \leq |\lambda_{n+1}| \int_{0}^{1} |k(t_1, s) - k(t_2, s)| |f(s, u_n(s), H[u_n])| \, ds \\
		& \leq c_\rho\max_{\Pi_\rho} |f| \int_{0}^{1} |k(t_1, s) - k(t_2, s)| \, ds \\
		& < c_\rho\max_{\Pi_\rho} |f| \varepsilon
	\end{align*}
	for each $n\in\N_0$.\\
	At the same time, by continuity of $u_0$, for every $\varepsilon>0$ there exists $\delta_0>0$ such that
	\[
	|u_0(t_1)-u_0(t_2)| < c_\rho\max_{\Pi_\rho} |f| \varepsilon
	\]
	for every $t_1,t_2\in[0,1]$ such that $|t_1-t_2|<\delta_0$. Therefore, we may take $\delta=\min\{\delta_0,\delta_1\}$ and this leads the equicontinuity.
	
	In summary, we have that $\{\lambda_n\}_n$ is a bounded sequence in $\R$ and $\{u_n\}_n$ is a equicontinuous and uniformly bounded sequence in $C[0,1]$. The Ascoli-Arzelà Theorem states that there exists a subsequence $\{u_{n_j}\}_j$ of $\{u_n\}_n$ uniformly converging to some $u^* \in C[0,1]$. Now, consider the subsequence $\{\lambda_{n_j}\}_j \subset \{\lambda_n\}_n$. This is a bounded sequence in $\R$, hence the Bolzano-Weierstrass Theorem yields the existence of another subsequence $\{\lambda_{n_{j_l}}\}_l$ that converges to some $\lambda^* \in \R$.\\
	Summarizing, we have two sequences $\{u_{n_{j_l}}\}_l$ and $\{\lambda_{n_{j_l}}\}_l$, where the first one uniformly converges to a continuous function $u^*$ and the second one converges to a real number $\lambda^*$. For the sake of notation, we will denote these two sequences again as $\{u_n\}_n$ and $\{\lambda_n\}_n$, respectively.\\
	We want to show that $(u^*,\lambda^*)$ is a couple satisfying
	\[
	\begin{cases}
		u^*(t) = \lambda^* Tu(t), &\text{for every } t \in [0,1],\\
		\norm{u^*} = \rho.
	\end{cases}
	\]
	First of all, by continuity of the norm,
	\[
	\norm{u^*}= \lim_{n\to\infty}\norm{u_n} = \lim_{n\to\infty}\rho = \rho.
	\]
	Moreover,
	\begin{align*}
		0 \leq \norm{u^*-\lambda^*Tu^*} &\leq \norm{u^*-u_n} + \norm{u_n-\lambda_nTu^*} + \norm{\lambda_nTu^*-\lambda^*Tu^*}\\
		&= \norm{u^*-u_n} + |\lambda_n|\norm{Tu_{n-1}-Tu^*} + |\lambda_n-\lambda^*|\norm{Tu^*}\\
		&\leq \norm{u^*-u_n} + c_\rho\norm{Tu_{n-1}-Tu^*} + |\lambda_n-\lambda^*|\norm{Tu^*} \\
		&\to 0 \quad\text{as } n \to \infty.
	\end{align*}
	This shows that the sequences that we have chosen converges, up to a subsequence, to the desired couple $(u^*,\lambda^*)$. 
	\section{Pseudocodes and MATLAB implementations}
	In this Section, we propose two pseudocodes. The first represents an approximation of the eigenvalues $\lambda_\rho$ together with the theoretical estimates presented in Theorem~\ref{HIEresult}, while the second implements the fixed-point iteration introduced in the Section~3, it is as a subroutine of the first and let us evaluate the approximation of the eigenpairs $(u_\rho^+,\lambda_\rho^+)$ and $(u_\rho^-,\lambda_\rho^-)$ for a fixed $\rho$. The pseudocodes are followed by a MATLAB implementation of them. Both the first pseudocode and the MATLAB implementation are in the case $\epsilon=-1$, but it is enough to change the Green's matrix, according to \eqref{greenneumannplus}, and the integral functions $\mathcal{A}$, $\mathcal{B}$, $\mathcal{C}$ and $\mathcal{D}$ for the case $\epsilon=1$, while the second pseudocode is applicable to both the cases $\epsilon=\pm 1$. Finally, at the end of this Section there are figures representing the eigenvalues approximations, their theoretical estimates and the approximations, for a fixed value of $\rho$, of the eigenfunctions, in both the cases $\epsilon=\pm 1$.
	
	\begin{algorithm}[H]
		\footnotesize
		\caption{\textsc{Approximation of the eigenvalues and theoretical estimates}}
		\begin{algorithmic}[1]
			\State $t \leftarrow \text{linspace}(0,1,N)$, \quad $\Delta s \leftarrow 1/(N-1)$
			\For{$i,j = 1,\ldots,N$}
			\State $G_{ij} \leftarrow \begin{cases}
				\dfrac{\cosh(1-t_j)\cosh(t_i)}{\sinh(1)} & \text{if } t_i \leq t_j \\[6pt]
				\dfrac{\cosh(1-t_i)\cosh(t_j)}{\sinh(1)} & \text{if } t_i > t_j
			\end{cases}$
			\EndFor
			
			\State Set $b \leftarrow 3\pi/2$, \quad $q \leftarrow \bigl(\sinh(1)(1+b^2)\bigr)^{-1}$
			\State Define the integral functions $\mathcal{A}(t), \mathcal{B}(t), \mathcal{C}(t), \mathcal{D}(t)$ on $[0,2/3]$ and $[2/3,1]$
			\State $\rho_1 \leftarrow \tfrac{1}{4}\max_{t \in [0,2/3]} \log(-\mathcal{A}(t)/\mathcal{B}(t))$
			\State $\rho_2 \leftarrow \tfrac{1}{4}\max_{t \in [2/3,1]} \log(-\mathcal{C}(t)/\mathcal{D}(t))$
			\State $\rho_0 \leftarrow \max\{\rho_1, \rho_2\}$
			
			\For{$\rho \in \text{linspace}(5\times10^{-3},\, m,\, 15)$, for any chosen $m \geq \rho_0$},
			\State $(\lambda_\rho^+, e_\rho^+) \leftarrow \textsc{FixedPointIteration}(\rho, +1, t, G, \Delta s)$
			\State $(\lambda_\rho^-, e_\rho^-) \leftarrow \textsc{FixedPointIteration}(\rho, -1, t, G, \Delta s)$
			\EndFor
			
			\For{$\rho \in \text{linspace}(0,\, \rho_0,\, 10^3)$}
			\State $\displaystyle\max_{[0,1]} \underline{F_\rho} \leftarrow \max\!\left(e^{-2\rho}\mathcal{A}(t)+e^{2\rho}\mathcal{B}(t),\; e^{-2\rho}\mathcal{C}(t)+e^{2\rho}\mathcal{D}(t)\right)$
			\State $a(\rho) \leftarrow \rho \,/\, \max_{[0,1]} \underline{F_\rho}$
			\EndFor
		\end{algorithmic}
	\end{algorithm}
	
	\begin{algorithm}[H]
		\footnotesize
		\caption{\textsc{FixedPointIteration}$(\rho,\, \sigma,\, t,\, G,\, \Delta s)$}
		\begin{algorithmic}[1]
			
			\State $u \leftarrow \rho \cdot \mathbf{1}$, \quad $\lambda \leftarrow \sigma$, where $\sigma=\operatorname{sign}(\lambda)$.
			
			\For{$k = 1, \ldots, \textit{max\_iter}$}
			\State $f \leftarrow \sin\!\left(\tfrac{3\pi}{2}t\right) \cdot e^u \,\Big/\, \int_0^1 e^{u(s)}\,ds$
			\State $w \leftarrow G \cdot f \cdot \Delta s$
			\State $\lambda_\text{new} \leftarrow \sigma\,\rho \,/\, \|w\|_\infty$
			\State $u_\text{new} \leftarrow \lambda_\text{new} \cdot w$
			\If{$\|u_\text{new} - u\|_\infty + |\lambda_\text{new} - \lambda| < \varepsilon$}
			\State $u \leftarrow u_\text{new}$, \quad $\lambda \leftarrow \lambda_\text{new}$
			\State \textbf{break}
			\EndIf
			\State $u \leftarrow u_\text{new}$, \quad $\lambda \leftarrow \lambda_\text{new}$
			\EndFor
			
			\State $f \leftarrow \sin\!\left(\tfrac{3\pi}{2}t\right) \cdot e^u \,\Big/\, \int_0^1 e^{u(s)}\,ds$
			\State $e \leftarrow \|u - \lambda\, G \cdot f \cdot \Delta s\|_\infty$
			\State \Return $(\lambda,\, e)$
		\end{algorithmic}
	\end{algorithm}
	
	Now, we present a MATLAB script that implements the fixed-point iteration in the last section to approximate the eigenvalues $\lambda_\rho^+>0$ and $\lambda_\rho^-$ for equation~\eqref{HIEneumann} in Example~\ref{Exampleminus} for every $\rho\in(0,\rho_0)$, where $\rho_0\approx0.2252$ is computed by the script itself.
	
\begin{figure}[H]
		\centering
		\begin{subfigure}[b]{0.95\textwidth}
			\centering
			\includegraphics[width=\linewidth]{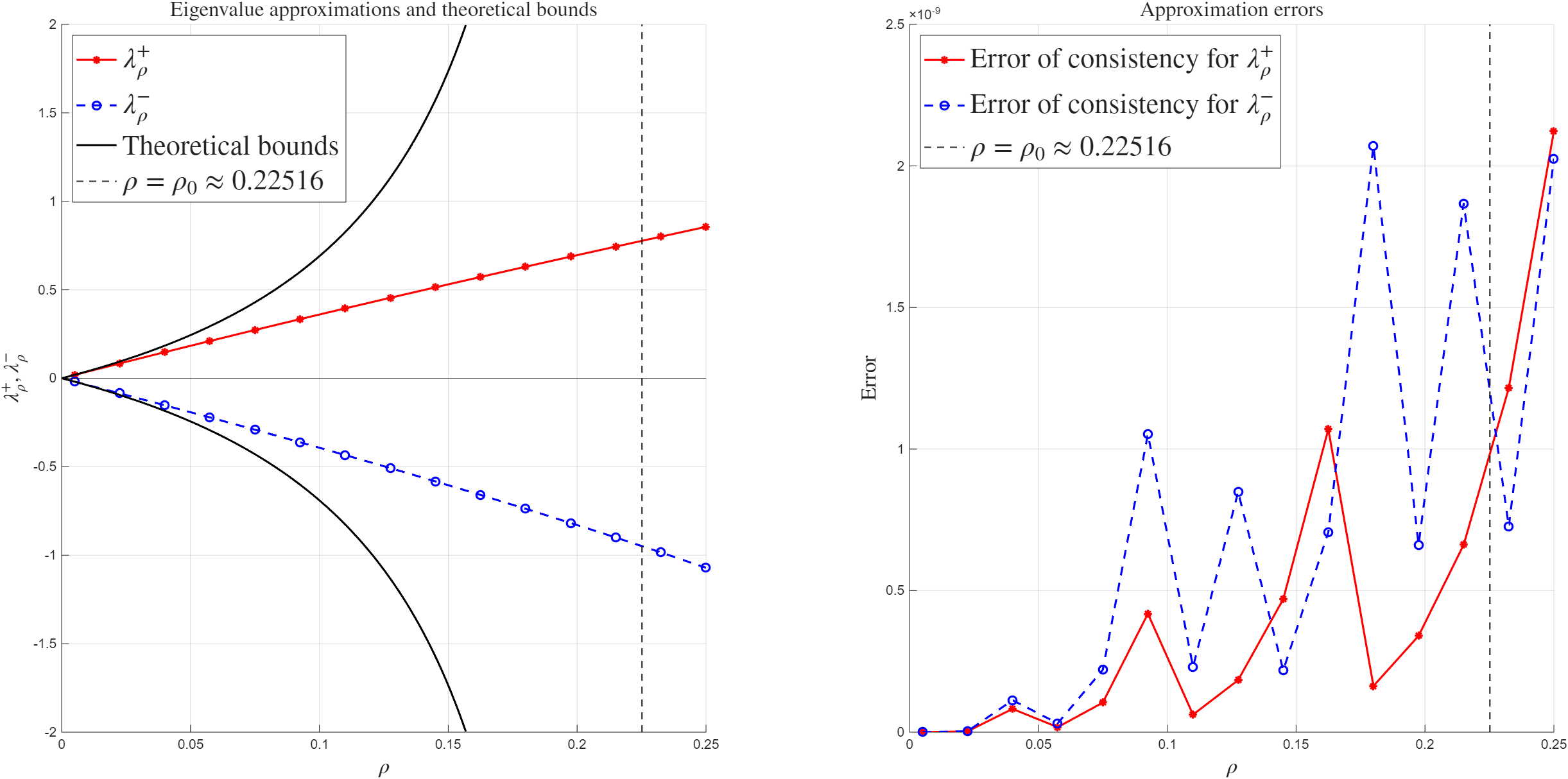}
			\caption{Case $\epsilon=-1$}
			\label{eigen_approx_minus}
		\end{subfigure}
		\hfill
		\begin{subfigure}[b]{0.95\textwidth}
			\centering
			\includegraphics[width=\linewidth]{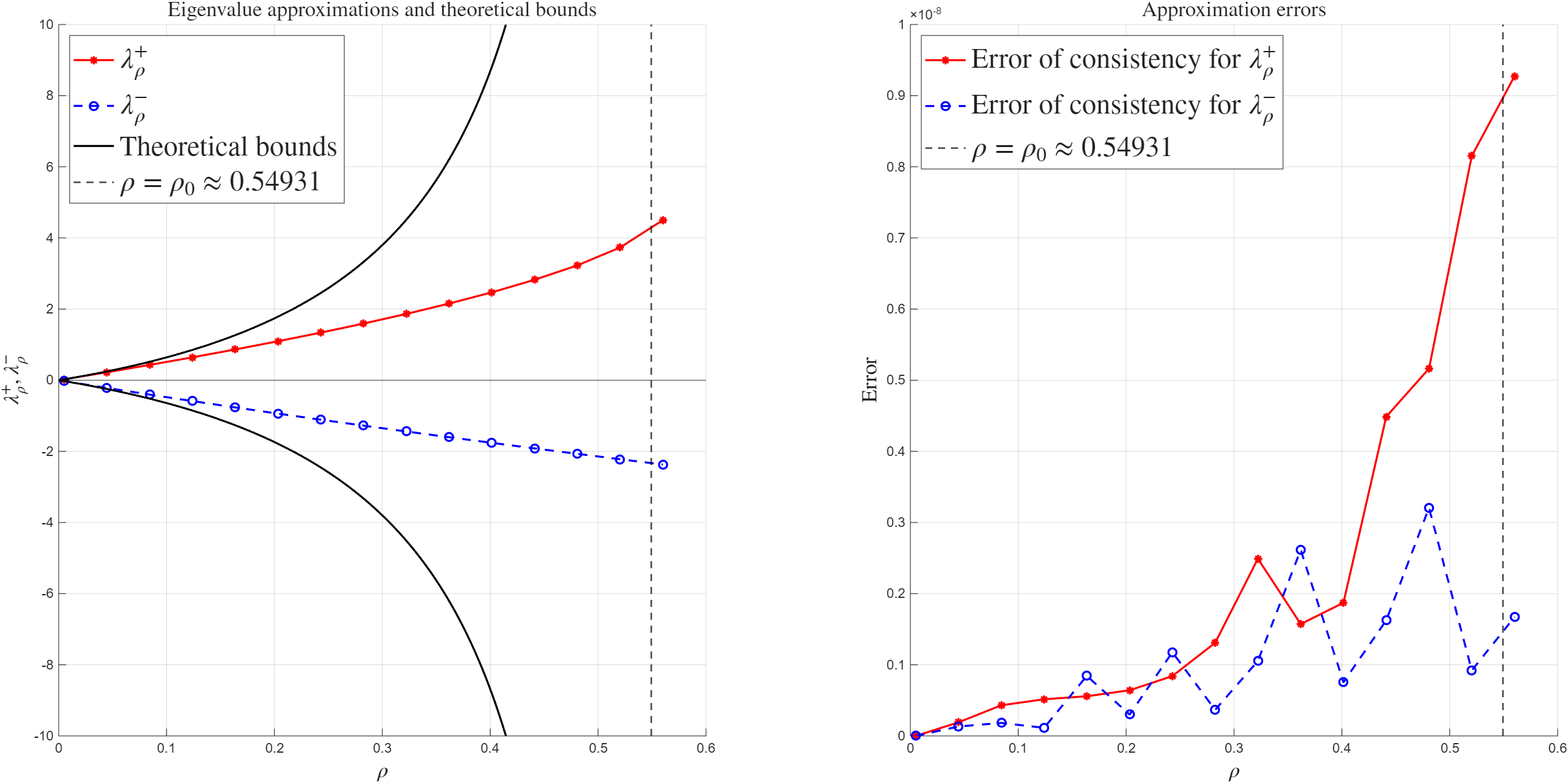}
			\caption{Case $\epsilon=1$}
			\label{eigen_approx_plus}
		\end{subfigure}
		\caption{The plots on the left illustrate the obtained approximations of the eigenvalues along with the theoretical bounds; the plots on the right illustrate the errors of consistency of the approximations.}
		\label{fig:eigen_approximations}
	\end{figure}

\begin{figure}[H]
		\centering
		\includegraphics[width=0.95\linewidth]{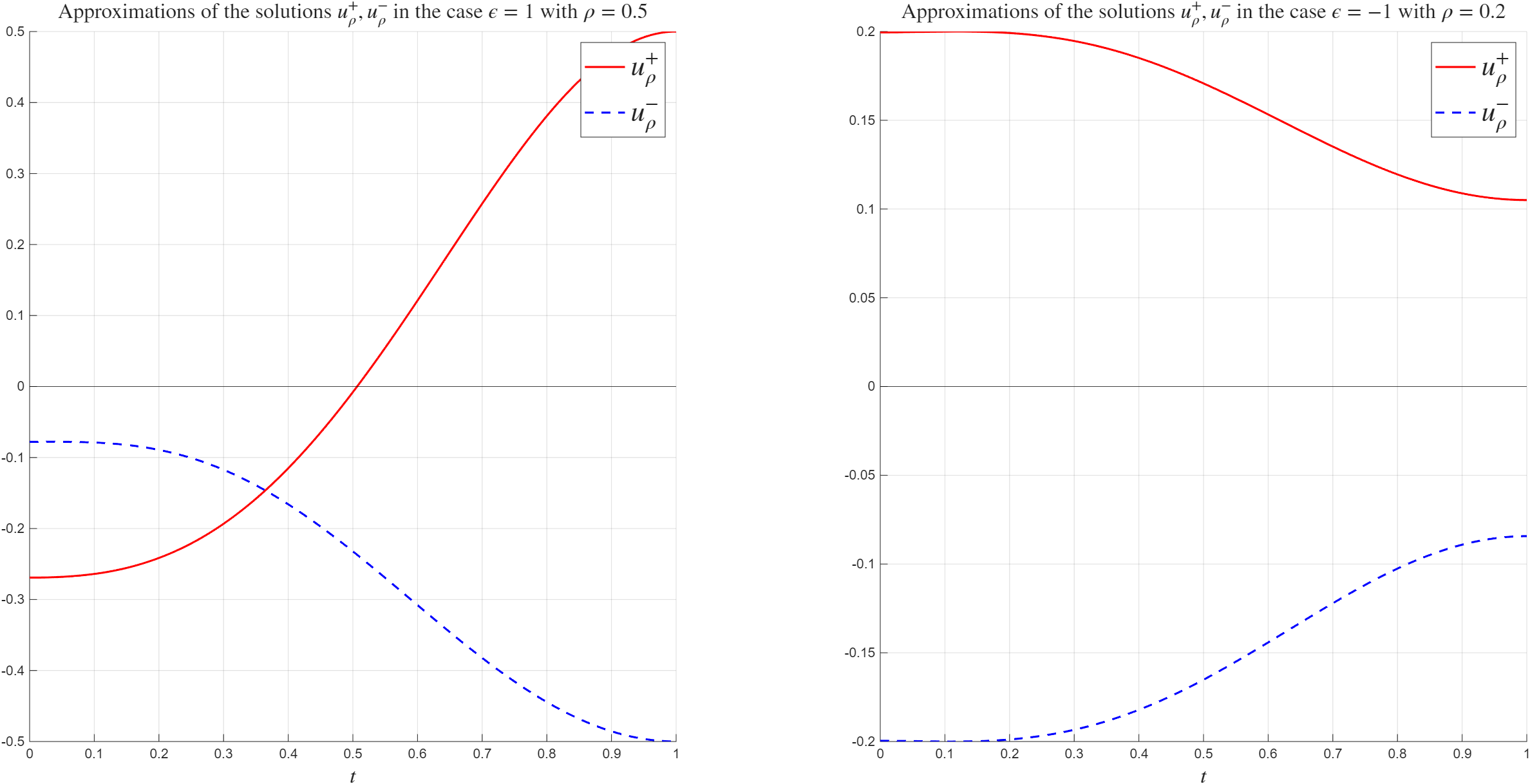}
		\caption{The plot on the left illustrates the approximation of the two eigenfunctions corresponding to $\rho=0.5$ in the case $\epsilon=1$; the plot on the right illustrates the approximation of the two eigenfunctions corresponding to $\rho=0.2$ in the case $\epsilon=-1$;}
		\label{fig:solution_approx_both_cases}
	\end{figure}

	\begin{lstlisting}[tabsize=1, keepspaces=true,numbers=left,style=Matlab-editor,breaklines=true, basicstyle=\footnotesize\linespread{1.1}\selectfont]
		%% Discretization
		t  = linspace(0,1,1e3)'; % Discrete points on [0,1]
		ds = 1/(length(t)-1);    % Size for each step
		
		%% Construction of the Green's matrix (homogeneous Neumann conditions)
		[T,S] = meshgrid(t,t);
		G     = (cosh(1-T).*cosh(S).*(S<=T)+...
		cosh(1-S).*cosh(T).*(S>T)) /sinh(1);
		
		%% Theoretical bounds and computation of maximal rho value
		% Parameters and auxiliary integral functions
		b  = 3*pi/2;
		q  = 1/(sinh(1)*(1+b^2));
		t1 = linspace(0,2/3,500);
		t2 = linspace(2/3,1,500);
		A  = @(t) q*(sin(b*t)*sinh(1)+b*(cosh(1-t)+cosh(1/3)*cosh(t)));
		B  = @(t) -q*b*cosh(t)*cosh(1/3);
		C  = @(t) q*b*cosh(1-t)*(cosh(2/3)+1);
		D  = @(t) q*(sin(b*t)*sinh(1)-b*cosh(1-t)*cosh(2/3));
		% Compute maximal rho value, i.e. rho0
		rho1 = max(log(-A(t1)./B(t1)))/4;
		rho2 = max(log(-C(t2)./D(t2)))/4;
		rho0 = max(rho1, rho2);
		% Theoretical curves
		rho_th    = linspace(0,.1568,1e3); % 0.1568 has been chosen for better plots
		F_low_max = zeros(size(rho_th));
		for i = 1:length(rho_th)
		r = rho_th(i);
		F_low = [exp(-2*r)*A(t1)+exp(2*r)*B(t1),...
		exp(-2*r)*C(t2)+exp(2*r)*D(t2)];
		F_low_max(i) = max(F_low);
		end
		a_rho = rho_th./F_low_max;
		
		%% Solver parameters
		max_iter = 1e3; % Maximum number of iterations
		tol = 1e-7;     % Tolerance for convergence
		rho = linspace(5e-3,.25,15); % Interval for rho, .25 has been chosen for better plots
		
		%% Solve both branches, l <-> lambda
		% Solve for positive eigenvalue (sgn = +1)
		[l_p,err_p] = arrayfun(@(r) solve(r,+1,t,G,ds,max_iter,tol), rho);
		% Solve for negative eigenvalue (sgn = -1)
		[l_n,err_n] = arrayfun(@(r) solve(r,-1,t,G,ds,max_iter,tol), rho);
		
		%% Plot #1: approximation of the eigenvalues, theoretical bounds, rho = rho0
		subplot(1,2,1)
		grid on, hold on
		plot(rho,l_p,'-*r',rho,l_n,'--ob')
		plot(rho_th,a_rho,'-k',rho_th,-a_rho,'-k')
		xline(rho0,'--k')
		yline(0, 'k')
		xlabel('\rho')
		ylabel('\lambda_\rho^+, \lambda_\rho^-')
		title('Eigenvalue approximations and theoretical bounds')
		legend('\lambda_\rho^+', '\lambda_\rho^-','Theoretical bounds', '',['\rho=\rho_0\approx', num2str(rho0)])
		
		%% Plot #2: residual errors, rho = rho0
		subplot(1,2,2)
		grid on, hold on
		plot(rho, err_p, '-*r', rho, err_n, '--ob')
		xline(rho0,'--k')
		xlabel('\rho')
		ylabel('Error')
		title('Approximation errors')
		legend('Error of consistency for \lambda_\rho^+', 'Error of consistency for \lambda_\rho^-',['\rho=\rho_0\approx', num2str(rho0)])
		
		%% Fixed-point iteration solver
		function [lambda,error] = solve(r,sgn,t,G,ds,max_iter,tol)
		u      = r*ones(size(t)); % Initial guess for u
		lambda = sgn;             % Initial guess for lambda
		for iter = 1:max_iter
		f          = sin(1.5*pi*t).*exp(u)/trapz(t,exp(u));
		w          = G*f*ds;
		lambda_new = sgn*r/norm(w,inf);
		u_new      = lambda_new*w;
		if norm(u_new-u,inf)+abs(lambda_new-lambda)<tol
		u      = u_new;
		lambda = lambda_new;
		break
		end
		u      = u_new;
		lambda = lambda_new;
		end
		f     = sin(1.5*pi*t).*exp(u)/trapz(t,exp(u));
		error = norm(u-lambda*G*f*ds,inf);
		end
	\end{lstlisting}
	The plots in Figure~\ref{fig:eigen_approximations} illustrate the approximations of the eigenvalues $\lambda_\rho^+,\lambda_\rho^-$ together with their theoretical estimates and the errors of consistency for the approximations, for both the cases $\epsilon=\pm 1$.
	The plots in Figure~\ref{fig:solution_approx_both_cases} show the approximations of the eigenfunctions $u_\rho^+,u_\rho^-$, for $\rho$ fixed, in both cases $\epsilon=\pm 1$.
	
	\section*{Acknowledgements}
	G.~A.~Veltri is a member of the ``Gruppo Nazionale per l'Analisi Matematica, la Probabilit\`a e le loro Applicazioni'' (GNAMPA) of the Istituto Nazionale di Alta Matematica (INdAM) and of the UMI Group TAA  ``Approximation Theory and Applications''.

\end{document}